\documentclass[11pt,a4paper]{article}
\usepackage{amsfonts}
\usepackage{mathrsfs}
\usepackage{amsmath}
\usepackage{dsfont}
\usepackage{graphicx,epstopdf}
\usepackage{amssymb}
\usepackage{theorem}
\usepackage{graphicx}
\usepackage[dvips]{epsfig}
\usepackage{latexsym}
\usepackage{exscale}
\usepackage[latin1]{inputenc}
\numberwithin{equation}{section}
\newtheorem{thm}{Theorem}[section]
\newtheorem{lemme}[thm]{Lemma}
\newtheorem{prop}[thm]{Proposition}

\newtheorem{defn}{Definition}[section]

\def\be#1\ee{\begin{equation}#1\end{equation}}
\newcommand{\bed}{\begin{displaymath}}
\newcommand{\eed}{\end{displaymath}}
\newcommand{\beq}{\begin{equation}}
\newcommand{\eeq}{\end{equation}}
\newcommand{\bea}{\begin{eqnarray}}
\newcommand{\eea}{\end{eqnarray}}
\newcommand{\beaa}{\begin{eqnarray*}}
\newcommand{\eeaa}{\end{eqnarray*}}
\newcommand{\bei}{\begin{itemize}}
\newcommand{\eei}{\end{itemize}}
\newcommand{\bee}{\begin{enumerate}}
\newcommand{\eee}{\end{enumerate}}

\def\Q{{\mathbb Q}}
\def\R{{\mathbb R}}
\def\E{{\mathbb E\,}}

\def\Z{{\mathbb Z}}
\def\N{{\mathbb Z}_{+}}

\def\e{\mathrm{e}}
\def\ud{\, \mathrm{d}}

\def\ga{{\gamma}}

\begin{document}
\title{H\"older regularity and series representation of a class of stochastic volatility models}
\author{Antoine Ayache and Qidi Peng\\
  Universit\'e Lille 1}
 \date{}
    \maketitle
\begin{abstract}
Let $\Phi:\R\rightarrow\R$ be an arbitrary continuously differentiable deterministic function such that $|\Phi|+|\Phi'|$ is bounded by a polynomial. 
In this article we consider the class of stochastic volatility models in which $\{Z(t)\}_{t\in [0,1]}$, the logarithm of the price process, is of the form $Z(t)=\int_{0}^t \Phi(X(s))\,dW(s)$, where $\{X(s)\}_{s\in[0,1]}$ denotes an arbitrary centered Gaussian 
process whose trajectories are, with probability $1$, H\"older continuous functions of an arbitrary order $\alpha\in (1/2,1]$, and where $\{W(s)\}_{s\in[0,1]}$ is a standard Brownian motion independent on $\{X(s)\}_{s\in [0,1]}$. 

First we show that the critical H\"older regularity of a typical trajectory of $\{Z(t)\}_{t\in[0,1]}$ is equal to $1/2$. Next we provide for  
such a trajectory an expression as a random series which converges at a geometric rate in any H\"older space of an arbitrary order $\ga<1/2$; this expression is obtained through the expansion of the random function $s\mapsto \Phi(X(s))$ on the Haar basis. Finally, thanks to it, we give an efficient iterative 
simulation method for $\{Z(t)\}_{t\in[0,1]}$.
\end{abstract}

\noindent{\small{\bf Running Title}:\ Series representations of stochastic volatility models}



\section{Introduction}
\label{sec:ch4-intro}

Stochastic volatility models are extensions of the
well-known Black and Scholes model. Hull and White \cite{HW} and
other authors working in the field of mathematical finance (see for instance \cite{Scott}
and \cite{MT}), introduced them starting in the mid-80s, in order to account for
volatility effects of exogenous arrivals of information on the price of an underlying asset. Later, for making such models to be
more realistic, it has been proposed (see for example \cite{CR1,CR2,GH1,GH,R}) 
to replace the Brownian motion governing the volatility in them, by a more
flexible stochastic process. That is why, in this article we consider the general class of stochastic volatility models, in which the stochastic process $\{Z(t)\}_{t\in[0,1]}$ denoting the logarithm of the price of an underlying asset, is defined for each time $t\in [0,1]$, as,
\begin{equation}
\label{eq0:ant-ch4}
Z(t)=\int_{0}^{t}\Phi(X(s))\ud W(s),
\end{equation}
where
\begin{itemize}
\item $\{X(s)\}_{s\in [0,1]}$ denotes an arbitrary real-valued centered Gaussian process whose trajectories are, with probability $1$,
H\"older continuous functions of an arbitrary order $\alpha\in (1/2,1]$ i.e. there is a positive and finite random variable $C$ such
that one has for almost all $\omega$ and for all $s_1,s_2\in [0,1]$,
\begin{equation}
\label{eq12:ant-ch4}
\big |X(s_1,\omega)-X(s_2,\omega)\big|\le C(\omega)|s_1-s_2|^{\alpha};
\end{equation}
it is worth noticing that, thanks to the Gaussianity of $\{X(s)\}_{s\in [0,1]}$, the random variable $C$ can be chosen in such a way that all its moments are finite.
\item $\Phi$ is an arbitrary deterministic real-valued continuously differentiable function on the real line, which satisfies as well as its derivative $\Phi'$, for all $x\in\R$,
\begin{equation}
\label{Cpol}
\big |\Phi (x)\big|+\big |\Phi' (x)\big|\le c\big(1+|x|\big)^{L},
\end{equation}
where $c>0$ and $L>0$ are two constants (non depending) on $x$.
\item
$\{W(s)\}_{s\in [0,1]}$ denotes a standard Brownian motion independent on the process $\{X(s)\}_{s\in [0,1]}$. Observe that thanks to these assumptions, the stochastic integral in (\ref{eq0:ant-ch4}) is well-defined, moreover, conditional on the $\sigma$-algebra, 
\begin{equation}
\label{art2:sigalG}
\mathcal {G}_X=\sigma\big(X(s),0\leq s\leq 1\big), 
\end{equation}
$\{Z(t)\}_{t\in[0,1]}$ is a centered Gaussian process with a covariance function given, for all $t_1,t_2\in[0,1]$, by,
\begin{equation}
\label{eq:lemme:2}
\E\big(Z(t_1)Z(t_2)|\mathcal{G}_X\big)=\int_0^{\min(t_1,t_2)}\big|\Phi(X(s))\big|^2\ud s.
\end{equation}
\end{itemize}
Observe that the models in (\ref{eq0:ant-ch4}) are generalizations of the multifractional stochastic volatility models introduced in \cite{Peng}, which in turn extend the 
fractional stochastic volatility models studied in \cite{GH1,GH}. Recall that a multifractional stochastic volatility model is defined 
through (\ref{eq0:ant-ch4}) in which $\{X(s)\}_{s\in [0,1]}$ denotes a multifractional Brownian motion \cite{BJR,PV}, in other words, roughly speaking, a fractional Brownian motion with a smooth time-varying Hurst parameter $H(s)$; a survey on the latter Gaussian process can be found in \cite{A} for instance, also we refer to \cite{Bi,BiP,BiPP} for some interesting applications of it in the field of finance. Recall that the main feature of a multifractional stochastic volatility model, is that the local roughness of the volatility process $\big\{\Phi(X(s)\big\}_{s\in [0,1]}$, 
in other words its pointwise H\"older exponent (see (\ref{thm:ch4-globalholferZ}) for a definition of this exponent), can be prescribed via the Hurst functional parameter $H(\cdot)$ and thus is allowed to change over time (that is from point to point).

The remaining of the paper is structured in the following way. In section~\ref{sec:ch4-holdreg}, we show that there exists a modification of $\{Z(t)\}_{t\in [0,1]}$, also denoted by
 $\{Z(t)\}_{t\in [0,1]}$, whose trajectories belong, with probability $1$, to any H\"older space $C^\gamma ([0,1])$ of an arbitrary order $\gamma<1/2$; on the other hand, under the assumption
 that $\Phi$ vanishes only on a Lebesgue negligible set, we show that the pointwise
 H\"older exponent of $\{Z(t)\}_{t\in [0,1]}$, at any point $t_0\in [0,1]$ such that $\E|X(t_0)|^2>0$, is almost surely equal to $1/2$; observe that, under the same 
 assumption on $\Phi$, the latter result implies that, with probability~$1$, the trajectories of $\{Z(t)\}_{t\in [0,1]}$ fail to belong to $C^\gamma ([0,1])$ when $\gamma>1/2$ and $\{X(s)\}_{s\in [0,1]}$ is not an almost surely vanishing process (i.e. we do not have for all $s\in [0,1]$, almost surely, $X(s)=0$).
 In section~\ref{sec:ch4-serrep}, by expanding the random function $s\mapsto \Phi(X(s))$ on the Haar basis of $L^2([0,1])$, we introduce a random series representation of
 $\{Z(t)\}_{t\in [0,1]}$, for which the convergence holds at a geometric rate, almost surely in $C^\gamma ([0,1])$, for all
 $\gamma<1/2$. Finally, thanks to the latter nice representation of $\{Z(t)\}_{t\in [0,1]}$, we give in section~\ref{sec:ch4-simulation} an iterative algorithm which allows to efficiently simulate this process.

 \section{H\"older regularity of the log price process $\{Z(t)\}_{t\in [0,1]}$}
 \label{sec:ch4-holdreg}

 Let us first recall the definition of a H\"older space of order $\gamma\in [0,1]$.

 \begin{defn}
\label{Holder space} For any $\gamma\in[0,1]$, the H\"older space
$C^{\gamma}([0,1])$ is defined as the Banach space of the continuous real-valued functions $u$ which satisfy,
$$
\sup_{0\le t_1 <t_2\le 1}\frac{|u(t_1)-u(t_2)|}{|t_1-t_2|^{\gamma}}<+\infty.
$$
It is equipped with the norm,
\begin{equation}
\label{eq:holdnorm-ch4}
\|u\|_{C^{\gamma}([0,1])}=\|u\|_{\infty}+\sup_{0\le t_1 <t_2\le 1}\frac{|u(t_1)-u(t_2)|}{|t_1-t_2|^{\gamma}},
\end{equation}
where
$
\|u\|_{\infty}=\sup_{t\in[0,1]}|u(t)|.
$
\end{defn}

Observe that $C^0 ([0,1])$ reduces to $C([0,1])$ the usual Banach space of real-valued continuous functions over $[0,1]$. The main goal of this section is to prove the following two theorems.

\begin{thm}
\label{thm:ch4-globalholferZ}
Let $\{Z(t)\}_{t\in [0,1]}$ be the stochastic process defined in (\ref{eq0:ant-ch4}), then there exists a modification of $\{Z(t)\}_{t\in [0,1]}$, also denoted by $\{Z(t)\}_{t\in [0,1]}$, such that, with probability $1$, for all $\gamma <1/2$, the trajectories of $\{Z(t)\}_{t\in [0,1]}$ belong to the H\"older space $C^\gamma ([0,1])$.
\end{thm}
Recall that when a trajectory $t\mapsto Y(t,\omega)$ of an arbitrary real-valued stochastic process $\{Y(t)\}_{t\in [0,1]}$
is a continuous and non-differentiable function at a point $t_0\in [0,1]$, then $\rho_Y (t_0,\omega)$, the pointwise H\"older exponent of $t\mapsto Y(t,\omega)$ at $t_0$, is defined as,
\begin{equation}
\label{eq13:ant-ch4}
\rho_Y (t_0,\omega)=\sup\left\{\rho\in [0,1]\,:\,\, \limsup_{t_0+h\in [0,1],\,\,\,h\rightarrow 0} \frac{Y(t_0+h,\omega)-Y(t_0,\omega)}{|h|^\rho}=0\right\}.
\end{equation}

\begin{thm}
\label{thm:1} Let $\{Z(t)\}_{t\in[0,1]}$ be the modification of the process $\{Z(t)\}_{t\in [0,1]}$ introduced in Theorem~\ref{thm:ch4-globalholferZ} and let $\rho_Z$ be the corresponding pointwise H\"older exponent; then, under the assumption that $\Phi$ vanishes only on a Lebesgue negligible set, one has,
\begin{equation}
\label{pHeZ1}
\mathbb{P}\big(\rho_{Z}(t_0)=1/2\big)=1,\, \mbox{$\forall$ $t_0\in[0,1]$ s.t. $\E|X(t_0)|^2>0$.}
\end{equation}
\end{thm}

Observe that a straightforward consequence of Theorem~\ref{thm:ch4-globalholferZ} is that, one has for all $t_0\in [0,1]$,
\begin{equation}
\label{eq114:ant-ch4}
\mathbb{P}\big(\rho_{Z}(t_0)\ge 1/2\big)=1. 
\end{equation}
Therefore, for obtaining Theorem~\ref{thm:1}, it is sufficient to show that Theorem~\ref{thm:ch4-globalholferZ} holds and that, under the assumption that $\Phi$ vanishes only on a Lebesgue negligible set, one has,
\begin{equation}
\label{eq114:ant-ch4}
\mathbb{P}\big(\rho_{Z}(t_0)\le 1/2\big)=1,\, \mbox{$\forall$ $t_0\in[0,1]$ s.t. $\E|X(t_0)|^2>0$.}
\end{equation}
The proof of Theorem~\ref{thm:ch4-globalholferZ} mainly relies on the following two lemmas. Observe that the first one whose proof can be found in \cite{KarShr} for instance, is a refined version of the usual Kolmogorov criterion allowing to show that a stochastic process has a modification with almost surely continuous trajectories; also observe that the second one, which can easily be proved, is a standard result on centered real-valued Gaussian random variables.

\begin{lemme}[Kolmogorov-Centsov criterion]
\label{theorem:Kol}
Assume that the positive real number $T$ is arbitrary and fixed. Let $\{M(t)\}_{t\in[0,T]}$ be
an arbitrary stochastic process such that for all $(s_1,s_2)\in[0,T]^2$, one has,
$$
\mathbb{E}|M(s_1)-M(s_2)|^{\tau}\leq c|s_1-s_2|^{1+\beta},
$$
where $c$, $\tau$ and $\beta$ are three positive constants. Then there exists a
modification of $\{M(t)\}_{t\in[0,T]}$ denoted by $\{\tilde{M}(t)\}_{t\in[0,T]}$, whose trajectories are
with probability~$1$, H\"older functions of any arbitrary order $\gamma\in[0,\beta/\tau)$; in other words,
$$
\mathbb{P}\left(\sup_{0\le s_1<s_2\le T}\frac{|\tilde{M}(s_1)-\tilde{M}(s_2)|}{|s_1-s_2|^{\gamma}}<+\infty\right)=1.
$$
\end{lemme}

\begin{lemme}[equivalence of Gaussian moments] Let $G$ be an arbitrary centered real-valued Gaussian random variable. Then for 
all positive real number $\tau$, one has,
\begin{equation}
\label{eq:gausmo}
\E |G|^\tau =c_\tau \big(\E |G|^2\big)^{\tau/2},
\end{equation}
where $c_\tau$ denotes the positive constant, only depending on $\tau$, defined as,
$$
c_\tau=\frac{2^{\tau/2}\Gamma\left(\frac{\tau+1}{2}\right)}{\Gamma\left(\frac{1}{2}\right)},
$$
$\Gamma$ being the usual "Gamma function", defined as $\Gamma(u)=\int_{0}^{+\infty} x^{u-1} e^{-x}\,dx$, for all $u\in (0,+\infty)$.
\end{lemme}

{\bf Proof of Theorem~\ref{thm:ch4-globalholferZ}:} Let $\gamma\in (0,1/2)$ be arbitrary and fixed, we assume that the positive real number $\tau$ has been chosen 
in such a way that
\beq
\label{eqnew:tauga}
\gamma<\frac{1}{2}-\frac{1}{\tau}.
\eeq
Recall that the $\sigma$-algebra $\mathcal {G}_X$ has been defined in (\ref{art2:sigalG}) and that 
conditionally on it, the stochastic process $\{Z(t)\}_{t\in [0,1]}$ has a centered Gaussian distribution with a covariance function given by (\ref{eq:lemme:2}).
Therefore, one has, for any real number $\tau>0$ and for each $t_1,t_2\in[0,1]$,
\begin{eqnarray}
\label{EE1}
\E|Z(t_1)-Z(t_2)|^{\tau}&=&\E\Big\{\E\big(|Z(t_1)-Z(t_2)|^{\tau}|\mathcal{G}_X\big)\Big\}\nonumber\\
&=& c_1\E\Big\{\Big(\E\big(|Z(t_1)-Z(t_2)|^{2}|\mathcal{G}_X\big)\Big)^{\tau/2}\Big\},
\end{eqnarray}
where $c_1>0$ is the constant $c_{\tau}$ in (\ref{eq:gausmo}), moreover, 
\begin{eqnarray}
\label{EE2}
&&\E\big(|Z(t_1)-Z(t_2)|^{2}|\mathcal{G}_X\big)\nonumber\\
&&=\E\big(|Z(t_1)|^{2}|\mathcal{G}_X\big)+\E\big(|Z(t_2)|^{2}|\mathcal{G}_X\big)-2\E\big(Z(t_1)Z(t_2)|\mathcal{G}_X\big)\nonumber\\
&&=\int_0^{t_1}\big|\Phi(X(s))\big|^2\ud s+\int_0^{t_2}\big|\Phi(X(s))\big|^2\ud s-2\int_0^{\min(t_1,t_2)}\big|\Phi(X(s))\big|^2\ud s\nonumber\\
&&=\int_{\min(t_1,t_2)}^{\max(t_1,t_2)}\big|\Phi(X(s))\big|^2\ud s.
\end{eqnarray}
Next, in order to bound $\E\big(|Z(t_1)-Z(t_2)|^{2}|\mathcal{G}_X\big)$, one sets 
\begin{equation}
\label{eqnouv-normX}
\|X\|_{\infty}=\sup_{t\in[0,1]}\big|X(s)\big|;
\end{equation}
observe that $\|X\|_{\infty}$ is an almost surely finite random variable, since the trajectories of $\{X(s)\}_{s\in [0,1]}$ are, with probability~$1$ continuous functions.
The Gaussianity of the latter process as well as the almost sure finiteness of $\|X\|_{\infty}$, imply that (see \cite{LedTal} for instance), all the moments of this random variable are finite as well, namely, for each real number $p>0$, one has,
\begin{equation}
\label{eq:ant-LT}
\E\big (\|X\|_{\infty}^p\big ) <\infty.
\end{equation}
Next, one sets,
\begin{equation}
\label{eq6:ant-ch4}
C_1=\sup_{x\in[-\|X\|_{\infty},\|X\|_{\infty}]}\big|\Phi(x)\big |,
\end{equation}
observe that $C_1$ is an almost surely finite random variable since $\Phi$ is a continuous function, moreover, in view of (\ref{Cpol}) and (\ref{eq:ant-LT}) all its moments
are finite, namely, for each real number $q>0$, one has,
\beq
\label{eq6:ant-ch4ter}
\E\big (C_{1}^q\big)<\infty.
\eeq
Next, combining (\ref{EE2}) with (\ref{eqnouv-normX}) and (\ref{eq6:ant-ch4}), one gets that,
\begin{equation}
\label{EE3}
\E\big(|Z(t_1)-Z(t_2)|^{2}|\mathcal{G}_X\big)\leq C_1^2|t_1-t_2|,
\end{equation}
therefore (\ref{EE1}) entails that
\begin{equation}
\label{EE4}
\E|Z(t_1)-Z(t_2)|^{\tau}\le c_2|t_1-t_2|^{\tau/2},
\end{equation}
where $c_2= c_1\big (\E|C_1|^{\tau}\big)<\infty$; observe that the finiteness of the latter constant is due to (\ref{eq6:ant-ch4ter}). Finally, putting together (\ref{EE4}), (\ref{eqnew:tauga}), and Lemma~\ref{theorem:Kol} in which one takes $T=1$ and $\beta=\tau/2-1$, one obtains the theorem. $\square$\\

From now on, our goal is to show that Relation (\ref{eq114:ant-ch4}) holds, to this end one needs the following lemma.
\begin{lemme}
\label{lemme:3} Assume that $t_0\in [0,1]$ is arbitrary and fixed. Let $(h_n)_n$ is an arbitrary sequence of non
vanishing real numbers which converges to $0$ and satisfies, for every $n$, $t_0+h_n\in [0,1]$. Then, 
\begin{itemize}
\item[(i)] for all fixed $\epsilon>0$, conditional on the $\sigma$-algebra $\mathcal
{G}_X$ defined in (\ref{art2:sigalG}), the random variable
$$
R_{t_0,h_{n}}=\frac{Z(t_0+h_n)-Z(t_0)}{|h_n|^{1/2+\epsilon}}
$$
has a centered Gaussian distribution with a variance given by,
$$
\mbox{var}\big(R_{t_0,h_{n}}|\mathcal
{G}_X\big)=|h_n|^{-1-2\epsilon}\int_{\min(t_0,t_0+h_n)}^{\max(t_0,t_0+h_n)}\big|\Phi(X(s))\big |^2\ud s;
$$
\item[(ii)] moreover, under the assumptions that $\Phi$ does not vanish except on a Lebesgue negligible set and $E|X(t_0)|^2>0$, one has, almost surely,
$$
\mbox{var}\big(R_{t_0,h_{n}}|\mathcal
{G}_X\big)\xrightarrow[n\rightarrow+\infty]{a.s.} +\infty.
$$
\end{itemize}
\end{lemme}
\textbf{Proof of Lemma \ref{lemme:3}:} Part~$(i)$ follows from (\ref{EE2}) and the fact that conditionally on $\mathcal
{G}_X$, the stochastic process $\{Z(t)\}_{t\in [0,1]}$ has a centered Gaussian distribution. Let us now prove that Part~$(ii)$ holds. 
The Mean Value Theorem, implies that, that there is a real number $\tilde{s}_n$ which belongs to the interval 
$\big(\min(t_0,t_0+h_n),\max(t_0,t_0+h_n)\big)$ and satisfies,
\begin{equation}
\label{assam1}
\int_{\min(t_0,t_0+h_n)}^{\max(t_0,t_0+h_n)}\big|\Phi(X(s))\big|^2\ud
s=|h_n|\big|\Phi(X(\tilde{s}_n))\big|^2;
\end{equation}
thus, using the almost sure continuity of the random function $s\mapsto \big|\Phi(X(s))\big|^2$, one gets that,
\begin{equation}
\label{conv:var}
|h_n|^{-1}\int_{\min(t,t+h_n)}^{\max(t,t+h_n)}\big|\Phi(X(s))\big|^2\ud
s\xrightarrow[n\rightarrow+\infty]{a.s.}\big|\Phi(X(t_0))\big|^2.
\end{equation}
Next, observe that the assumptions that $\Phi$ vanishes only on Lebesgue negligible set and that $\E|X(t_0)|^2>0$, imply that,
almost surely,
\begin{equation}
\label{assam2}
\big|\Phi(X(t_0))\big|^2 \stackrel{a.s.}{>}0.
\end{equation}
Finally it follows from (\ref{conv:var}) and (\ref{assam2}) that
$$
|h_n|^{-1-2\epsilon}\int_{\min(t_0,t_0+h_n)}^{\max(t_0,t_0+h_n)}\big|\Phi(X(s))\big|^2\ud
s\xrightarrow[n\rightarrow+\infty]{a.s.}+\infty.
$$
 $\square$\\
 
 Now, we are in position to prove Relation (\ref{eq114:ant-ch4}).\\
 {\bf Proof of Relation (\ref{eq114:ant-ch4}):} Our proof is inspired by that of Proposition~2.4 in \cite{AL:2000}. It consists in showing that for all fixed $\epsilon>0$, there exists a  sequence $(r_k)_{k\in\mathbb{N}}$ of non vanishing real numbers converging to $0$ which satisfies:
\begin{equation}
\label{diff:Z}
\mbox{$t_0+r_k\in [0,1]$ for all $k$, and }\, \frac{|Z(t_0+r_k)-Z(t_0)|}{|r_k|^{1/2+\epsilon}}\xrightarrow[k\rightarrow+\infty]{a.s.}+\infty.
\end{equation}
To this end, it is sufficient to prove that there exists a sequence $(h_n)_n$ of non vanishing real numbers converging to $0$ which satisfies:
\begin{equation}
\label{diff:Z1}
\mbox{$t_0+h_n\in [0,1]$ for all $n$, and }\, \frac{1}{R_{t_0,h_n}}=\frac{|h_n|^{1/2+\epsilon}}{|Z(t_0+h_n)-Z(t_0)|}\xrightarrow[n\rightarrow+\infty]{\mathbb{P}}0,
\end{equation}
where the notation "$\stackrel{\mathbb{P}}{\longrightarrow} 0$" means that the convergence to $0$ holds in probability. Indeed, assuming that (\ref{diff:Z1}) is satisfied, then one can extract from $(h_n)_n$ a subsequence denoted by $(r_k)_{k\in\N}$ such that one has (\ref{diff:Z}).
Let us now prove (\ref{diff:Z1}). Denote by $(h_n)_n$ an arbitrary sequence of non vanishing real numbers converging to $0$ and such that $t_0+h_n\in [0,1]$ for all $n$.
Observe that, one has, for all real number $\eta > 0$,
\begin{eqnarray}
\label{PE}
\mathbb{P}\left(\frac{1}{R_{t_0,h_n}}\leq\eta\right) =\mathbb{E}\Big(\mathds{1}_{\big(\frac{1}{R_{t_0,h_n}}\leq\eta\big)}\Big)&=&\mathbb{E}\bigg(\mathbb{E}\Big(\mathds{1}_{\big(\frac{1}{R_{t_0,h_n}}\leq\eta\big)}\Big|\mathcal{G}_X\Big)\bigg)\nonumber\\
&=&\mathbb{E}\bigg(\mathbb{E}\Big(\mathds{1}_{\big (R_{t_0,h_n}\geq 1/\eta\big)}\Big|\mathcal{G}_X\Big)\bigg).\nonumber\\
\end{eqnarray}
Moreover, Lemma~\ref{lemme:3} entails that, almost surely,
\begin{equation}
\label{eq3:ant-ch4}
\mathbb{E}\Big(\mathds{1}_{\big (R_{t_0,h_n}\geq 1/\eta\big)}\Big|\mathcal{G}_X\Big)\stackrel{a.s.}{=}\sqrt{\frac{2}{\pi}}\int_{\eta^{-1}\left (\mbox{var}\left(R_{t_0,h_{n}}|\mathcal {G}_X\right)\right)^{-1/2}}^{+\infty} \exp(-x^2/2)\ud x
\end{equation}
and
\begin{equation}
\label{eq2:ant-ch4}
\mbox{var}\left(R_{t_0,h_{n}}|\mathcal {G}_X\right)\xrightarrow[n\rightarrow+\infty]{a.s.} +\infty;
\end{equation}
Thus, combining (\ref{eq3:ant-ch4}) with (\ref{eq2:ant-ch4}) one obtains that,
\begin{equation}
\label{eq4:ant-ch4}
\mathbb{E}\Big(\mathds{1}_{\big(R_{t_0,h_n}\geq 1/\eta\big)}\Big|\mathcal{G}_X\Big)\xrightarrow[n\rightarrow+\infty]{a.s.}\sqrt{\frac{2}{\pi}}\int_{0}^{+\infty} e^{-x^2/2}\ud x=1.
\end{equation}
Finally, in view of (\ref{PE}) and (\ref{eq4:ant-ch4}), using the dominated convergence theorem, it follows that (\ref{diff:Z1}) holds.
$\square$

\section{Random series representation of $\{Z(t)\}_{t\in [0,1]}$ via the Haar basis}
\label{sec:ch4-serrep}
In order to state the main result of this section, we need to introduce some notations.
\begin{itemize}
\item We denote by $L^2([0,1])$ the usual Lebesgue Hilbert of the square integrable real-valued deterministic functions over $[0,1]$.
\item The Haar orthonormal basis of $L^2([0,1])$ (see for example \cite{M1,D,W}), is the
sequence of the functions:
$$
\big\{\varphi_{0,0},\psi_{j,k}:\,j\in\N \mbox{ and } k\in\{0,\ldots, 2^{j}-1\}\big\},
$$
defined as
\begin{equation}
\label{eq16:ant-ch4}
\varphi_{0,0}=\mathds{1}_{[0,1)}
\end{equation}
and 
\begin{equation}
\label{eq17:ant-ch4}
\psi_{j,k}=2^{j/2}(\mathds{1}_{[\frac{k}{2^j},\frac{2k+1}{2^{j+1}})}-\mathds{1}_{[\frac{2k+1}{2^{j+1}},\frac{k+1}{2^j})}),
\end{equation}
where $\mathds{1}_{S}$ is the indicator function of an arbitrary set $S$.
\item We denote by $(\Omega,\mathcal{F},\mathbb{P})$ the underlying probability space, that is
the probability space on which the processes $\{X(s)\}_{s\in [0,1]}$, $\{W(s)\}_{s\in [0,1]}$ and
$\{Z(t)\}_{t\in [0,1]}$ are defined; moreover, for the sake of simplicity, we assume that this space has been chosen in such a way that: {\em for all}
$\omega\in\Omega$, Relation (\ref{eq12:ant-ch4}) holds and $Z(\cdot,\omega)\in C^{\gamma}([0,1])$ for each $\gamma\in [0,1/2)$. Also, we denote by 
$$
\big\{\delta_{0,0},\lambda_{j,k}:\,j\in\N \mbox{ and } k\in\{0,\ldots, 2^{j}-1\}\big\},
$$ 
the sequence of standard independent Gaussian random variables defined, on this space as,
\begin{equation}
\label{eq13:ant-ch4}
\delta_{0,0}=\int_0^1\varphi_{0,0}(s)\ud W(s)=W(1)
\end{equation}
 and 
\begin{equation}
\label{eq14:ant-ch4}
\lambda_{j,k}=\int_0^1 \psi_{j,k}(s)\ud W(s)=-2^{j/2}\Big(W(\frac{k+1}{2^{j}})-2W(\frac{2k+1}{2^{j+1}})+W(\frac{k}{2^{j}})\Big).
\end{equation}
Observe that, similarly to Lemma~2 in \cite{AT}, one can show that there exist $\Omega^*\subseteq \Omega$
an event of probability $1$, and a positive random variable $C_{*}$ of finite moment of any order,
such that one has for all $\omega\in\Omega^*$, all $j\in\N$ and all $k\in\{0,\ldots, 2^j-1\}$,
\begin{equation}
\label{lambda} |\lambda_{j,k}(\omega)|\leq C_* (\omega)\sqrt{1+j};
\end{equation}
we note in passing that the proof of this important inequality, mainly relies on Borel-Cantelli Lemma.
\item $\{K(t,s)\}_{(t,s)\in [0,1]^2}$ is the stochastic
field defined for all $(t,s,\omega)\in [0,1]^2 \times \Omega$ as,
\begin{equation}
\label{eq5:ant-ch4}
K(t,s,\omega):=\Phi(X(s,\omega))\mathds{1}_{[0,t]}(s);
\end{equation}
thus, the process $\{Z(t)\}_{t\in [0,1]}$ defined in (\ref{eq0:ant-ch4}) can be expressed as,
\begin{equation}
\label{Z1} Z(t)=\int_0^1K(t,s)\ud W(s).
\end{equation}
Observe that (\ref{eq5:ant-ch4}), (\ref{eqnouv-normX}), (\ref{eq6:ant-ch4}), imply that for all $\omega\in\Omega$,
\begin{equation}
\label{eq11:ant-ch4bis}
\sup_{(t,s)\in [0,1]^2} \big|K(t,s,\omega)\big |=\sup_{s\in [0,1]} \big |\Phi(X(s,\omega))\big|\leq C_1(\omega)<\infty,
\end{equation}
where the last inequality results from the continuity of the functions $s\mapsto X(s,\omega)$ and $x\mapsto \Phi(x)$.
\item We denote by $\{b_{0,0}(t)\}_{t\in [0,1]}$ the stochastic process defined for all
$(t,\omega)\in [0,1]\times\Omega$, as,
\begin{equation}
\label{eq7:ant-ch4}
b_{0,0}(t,\omega)=\int_0^1K(t,s,\omega)\varphi_{0,0}(s)\ud s;
\end{equation}
moreover, for all $j\in \N$ and $k\in\{0,\dots, 2^j-1\}$, we denote by $\{a_{j,k}(t)\}_{t\in [0,1]}$ the stochastic process defined for all
$(t,\omega)\in [0,1]\times\Omega$, as,
\begin{equation}
\label{eq8:ant-ch4}
a_{j,k}(t,\omega)=\int_0^1 K(t,s,\omega)\psi_{j,k}(s)\ud s.
\end{equation}
Observe that it follows from (\ref{eq7:ant-ch4}), (\ref{eq8:ant-ch4}),
(\ref{eq11:ant-ch4bis}), (\ref{eq16:ant-ch4}) and (\ref{eq17:ant-ch4}), that for
all $\omega\in\Omega$ and for all $t_1,t_2\in [0,1]$,
\begin{equation}
\label{eq18:ant-ch4}
\big |b_{0,0}(t_1,\omega)-b_{0,0}(t_2,\omega)\big|\le C_1(\omega) |t_1-t_2|
\end{equation}
and for all $j\in \N$ and $k\in\{0,\dots, 2^j-1\}$,
\begin{equation}
\label{eq19:ant-ch4}
\big |a_{j,k}(t_1,\omega)-a_{j,k}(t_2,\omega)\big |\le C_1(\omega) 2^{j/2}|t_1-t_2|.
\end{equation}
\end{itemize}
Now we are in position to state the main result of this section:
\begin{thm}
\label{thm:z}
Let $\gamma\in [0,1/2)$ be arbitrary and fixed. For all $J\in\N$ and $(t,\omega)\in [0,1]\times\Omega$, one sets,
\begin{equation}
\label{eq15:ant-ch4}
Z_J(t,\omega)=b_{0,0}(t,\omega)\delta_{0,0}(\omega)+\sum_{j=0}^{J}\sum_{k=0}^{2^j-1}a_{j,k}(t,\omega)\lambda_{j,k}(\omega).
\end{equation}
In view of (\ref{eq18:ant-ch4}) and (\ref{eq19:ant-ch4}) the trajectories of the process $\{Z_J (t)\}_{t\in [0,1]}$ belong to the H\"older space $C^\gamma ([0,1])$ since they are in fact Lipschitz functions; moreover, there exist $\Omega_{2}^*$ an event of probability $1$ and a positive random variable $D$ of finite moment of any order, such that one has for
all $\omega\in\Omega_{2}^*$ and $J\in\N$,
\begin{equation}
\label{eq20:ant-ch4}
\big \|Z(\cdot,\omega)-Z_J (\cdot,\omega)\big\|_{C^\gamma ([0,1])}\le D(\omega) 2^{-J\min(1/2-\gamma,\alpha-1/2)}\,\sqrt{1+J},
\end{equation}
where $\|\cdot\|_{C^\gamma ([0,1])}$ is the usual norm on $C^\gamma ([0,1])$ (see Definition~\ref{Holder space}) and where $\alpha\in (1/2,1]$ has been introduced in
(\ref{eq12:ant-ch4}).
\end{thm}
In order to prove the latter theorem we need some preliminary results. The following lemma is a weak version of Theorem~\ref{thm:z}.
\begin{lemme}
\label{lem1:ant-ch4}
We use the same notations as in Theorem~\ref{thm:z}. Let $t\in [0,1]$ be arbitrary and fixed. When $J$ goes to $+\infty$, the random variable
$Z_J (t)$ converges to the random variable $Z(t)$ in the Hilbert space $L^2(\Omega)$, namely, one has,
\begin{equation}
\label{eq21:ant-ch4}
\lim_{J\rightarrow +\infty} \E\big (Z(t)-Z_J (t)\big)^2 =0.
\end{equation}
As a straightforward consequence, there exist $\Omega_{1,t}^{*}$ an event of probability~$1$ included in $\Omega^*$ (recall that $\Omega^*$ is the event of probability~$1$ on which (\ref{lambda}) holds) and a subsequence $n\mapsto J_n$ (a priori depending on $t$) such that for all $\omega\in \Omega_{1,t}^{*}$,
one has,
\begin{equation}
\label{eq22:ant-ch4}
\lim_{n\rightarrow +\infty} Z_{J_n}(t,\omega)=Z(t,\omega).
\end{equation}
\end{lemme}

\noindent {\bf Proof of Lemma~\ref{lem1:ant-ch4}:} Let $K$ be as in (\ref{eq5:ant-ch4}). Observe that in view of (\ref{eq11:ant-ch4bis}), for all fixed $(t,\omega)\in [0,1]\times \Omega$, the function $K(t,\cdot,\omega):s\mapsto K(t,s,\omega)$ belongs to $L^2([0,1];\ud s)$. By expanding the latter function on the Haar basis, one obtains that,
\begin{equation}
\label{dec:K}
K(t,\cdot,\omega)=b_{0,0}(t,\omega)\varphi_{0,0}(\cdot)+\sum_{j=0}^{+\infty}\sum_{k=0}^{2^j-1}a_{j,k}(t,\omega)\psi_{j,k}(\cdot),
\end{equation}
where the coefficients $b_{0,0}(t,\omega)$ and $a_{j,k}(t,\omega)$ have been defined respectively in (\ref{eq7:ant-ch4}) and
(\ref{eq8:ant-ch4}). A priori the series in (\ref{dec:K}) is convergent in the $L^2([0,1];\ud s)$ norm, namely,
\begin{equation}
\label{eq9:ant-ch4}
\lim_{J\rightarrow +\infty} \int_{0}^1 \big |K(t,s,\omega)-K_J (t,s,\omega)\big|^2 \ud s=0,
\end{equation}
where, for each $J\in\Z_+$, $K_J (t,\cdot,\omega)$ is the partial sum defined, for all $s\in [0,1]$, as,
\begin{equation}
\label{eq9bis:ant-ch4}
K_J (t,s,\omega)=b_{0,0}(t,\omega)\varphi_{0,0}(s)+\sum_{j=0}^{J}\sum_{k=0}^{2^j-1}a_{j,k}(t,\omega)\psi_{j,k}(s).
\end{equation}
Let us show that this series, is also convergent in the
$L^2\big([0,1]\times\Omega;\ud s \otimes \mathbb{P}\big)$ norm, that is,
\begin{equation}
\label{eq10:ant-ch4}
\lim_{J\rightarrow +\infty} \E\left (\int_{0}^1 \big |K(t,s)-K_J (t,s)\big|^2 \ud s\right)=0.
\end{equation}
In order to show that (\ref{eq10:ant-ch4}) holds, we will use the dominated convergence theorem. It follows from (\ref{dec:K}), (\ref{eq9bis:ant-ch4}), Parseval formula and 
(\ref{eq11:ant-ch4bis}), that for all $(\omega,J)\in\Omega\times\N$,
\begin{eqnarray}
\label{eq11:ant-ch4}
\int_{0}^1  \big |K(t,s,\omega)-K_J (t,s,\omega)\big|^2 \ud s \le \int_{0}^1 \big |K(t,s,\omega)\big|^2 \ud s \le C_{1}^2(\omega);\nonumber
\end{eqnarray}
thus, in view of (\ref{eq9:ant-ch4}) and (\ref{eq6:ant-ch4ter}), we are allowed to derive (\ref{eq10:ant-ch4}) by making use of the dominated convergence theorem.

Finally, observe that, (\ref{Z1}), (\ref{eq13:ant-ch4}), (\ref{eq14:ant-ch4}), (\ref{eq15:ant-ch4}), (\ref{eq9bis:ant-ch4}), (\ref{art2:sigalG}), and the isometry property of Wiener integral, entail that, one has almost surely, for all $J\in\Z_+$,
$$
\E\Big( \big (Z(t)-Z_J (t)\big)^2|\mathcal{G}_X\Big)=\int_{0}^1 \big |K(t,s)-K_J (t,s)\big|^2 \ud s,
$$
which in turn implies that,
\begin{eqnarray*}
\E\big (Z(t)-Z_J (t)\big)^2&=& E\left(\E\Big( \big (Z(t)-Z_J (t)\big)^2|\mathcal{G}_X\Big)\right)\\
&=& \E\left (\int_{0}^1 \big |K(t,s)-K_J (t,s)\big|^2 \ud s\right);
\end{eqnarray*}
therefore, using (\ref{eq10:ant-ch4}) one gets (\ref{eq21:ant-ch4}). $\square$\\

The following lemma provides sharp estimates for $a_{j,k}(t)$.
\begin{lemme}
\label{lem2:ant-ch4}
There is a positive random variable $A$ of finite moment of any order, such that for all $\omega\in\Omega$, $t\in [0,1]$,
$j\in\N$ and $k\in\{0,\ldots, 2^j-1\}$,
\begin{itemize}
\item[(i)] when
$t\ge (k+1)/2^j$, one has
\begin{equation}
\label{eq23:ant-ch4}
\big|a_{j,k}(t,\omega)\big|\le A(\omega) 2^{-j(\alpha+1/2)},
\end{equation}
where $\alpha$ has been introduced in (\ref{eq12:ant-ch4});
\item[(ii)] when $k/2^j< t < (k+1)/2^j$, one has
\begin{equation}
\label{eq24:ant-ch4}
\big|a_{j,k}(t,\omega)\big|\le A(\omega) 2^{j/2}\big(t-k/2^j\big);
\end{equation}
\item[(iii)] when $t\le k/2^j$, one has,
\begin{equation}
\label{eq25:ant-ch4}
a_{j,k}(t,\omega)=0
\end{equation}
\end{itemize}
\end{lemme}

\noindent \textbf{Proof of Lemma \ref{lem2:ant-ch4}:} First observe that
(\ref{eq25:ant-ch4}) easily results from (\ref{eq17:ant-ch4}), (\ref{eq5:ant-ch4}) and
(\ref{eq8:ant-ch4}). Let us now prove that (\ref{eq23:ant-ch4}) holds, so we assume that $t\ge
(k+1)/2^j$. For all $(y,\omega)\in [0,1]\times\Omega$, we set
\begin{equation}
\label{eq26:ant-ch4}
X_{\Phi}^{(-1)}(y,\omega)=\int_{0}^{y}\Phi(X(s,\omega))\ud s.
\end{equation}
Then, using (\ref{eq17:ant-ch4}), (\ref{eq5:ant-ch4}),
(\ref{eq8:ant-ch4}) and (\ref{eq26:ant-ch4}), it follows that $a_{j,k}(t,\omega)$ can be
expressed as an increment of order $2$ of the function $2^{j/2}X_{\Phi}^{(-1)}(\cdot,\omega)$:
namely, when $t\ge(k+1)/2^j$, one has,
\begin{eqnarray}
\label{eq27:ant-ch4}
&& a_{j,k}(t,\omega)\\
&& =2^{j/2}\Big (X_{\Phi}^{(-1)}\big(\frac{2k+1}{2^{j+1}},\omega\big)-X_{\Phi}^{(-1)}\big(\frac{k}{2^j},\omega\big)\Big)\nonumber\\
&& \hspace{5cm} -2^{j/2}\Big (X_{\Phi}^{(-1)}\big(\frac{k+1}{2^{j}},\omega\big)-X_{\Phi}^{(-1)}\big(\frac{2k+1}{2^{j+1}},\omega\big)\Big).\nonumber
\end{eqnarray}
Next, applying the Mean Value Theorem to the function
$$
y\mapsto
X_{\Phi}^{(-1)}\big(\frac{k+1}{2^{j}}-y,\omega\big)-X_{\Phi}^{(-1)}\big(\frac{2k+1}{2^{j+1}}-y,\omega\big),
$$
on the interval $[0,2^{-j-1}]$, it follows from (\ref{eq27:ant-ch4}) and (\ref{eq26:ant-ch4}), that there
exists $z\in (0,2^{-j-1})$ such that
\begin{equation}
\label{eq28:ant-ch4}
a_{j,k}(t,\omega)=2^{-j/2-1}\left(\Phi\Big(X\Big(\frac{2k+1}{2^{j+1}}-z,\omega\Big)\Big)-\Phi\Big(X\Big(\frac{k+1}{2^{j}}-z,\omega\Big)\Big)\right);
\end{equation}
then, applying the same theorem to the function $x\mapsto \Phi(x)$, on the
compact interval whose endpoints are
$X\Big(\frac{2k+1}{2^{j+1}}-z,\omega\Big)$ and
$X\Big(\frac{k+1}{2^{j}}-z,\omega\Big)$, one obtains, in view of (\ref{eq28:ant-ch4}) and (\ref{eqnouv-normX}), that,
\begin{equation}
\label{eq29:ant-ch4}
\big |a_{j,k}(t,\omega)\big |\le C_2 (\omega) 2^{-j/2-1} \left
  |X\Big(\frac{2k+1}{2^{j+1}}-z,\omega\Big)-X\Big(\frac{k+1}{2^{j}}-z,\omega\Big)\right|,
\end{equation}
where
$$
C_2(\omega):=\sup\Big\{\big |\Phi'(x)\big |\,:\,\, |x|\le \|X\|_{\infty}(\omega)\Big\}.
$$
Observe that the continuity of the functions $s\mapsto X(s,\omega)$ and $x\mapsto \Phi'(x)$, implies that $C_2$ is a finite random variable, moreover (\ref{Cpol}) and  (\ref{eq:ant-LT}) entail that all its moments are finite as well. Next, combining (\ref{eq29:ant-ch4}) with (\ref{eq12:ant-ch4}), one gets
(\ref{eq23:ant-ch4}). Let us now show that (\ref{eq24:ant-ch4}) holds. It follows from (\ref{eq17:ant-ch4}), (\ref{eq5:ant-ch4}) and (\ref{eq8:ant-ch4}) that
$$
\big|a_{j,k}(t,\omega)\big|\le 2^{j/2}\int_{k/2^j}^{t} \big|\Phi(X(s))\big|\ud s\le C_1(\omega) 2^{j/2}\big(t-k/2^j\big),
$$
where the random variable $C_1$ has been defined in (\ref{eq6:ant-ch4}); thus we obtain (\ref{eq24:ant-ch4}). $\square$

\begin{lemme}
\label{lem3:ant-ch4}
Let $\delta$ be an arbitrary fixed positive real number. There is a constant $c>0$, only depending on $\delta$, such that for all $x,y\in (0,1]$,
satisfying $x\le y$, one has
$$
x^{\delta}\,\sqrt{1+\log_2 (x^{-1})}\le c \,y^{\delta}\,\sqrt{1+\log_2 (y^{-1})}.
$$
\end{lemme}
\noindent{\bf Proof of Lemma~\ref{lem3:ant-ch4}:} The derivative over $(0,1]$ of the positive function $s\mapsto s^{\delta}\,\sqrt{1+\log_2 (s^{-1})}$ is equal
to
$$
s^{\delta-1}\left (1-\frac{\log s}{\log 2}\right )^{-1/2}\left (\delta-\delta\frac{\log s}{\log 2}-\frac{1}{2\log 2}\right ).
$$
Let us set $\overline{s}_\delta=2\e^{-(2\delta)^{-1}}$. When $\overline{s}_\delta\ge 1$, then the function $s\mapsto s^{\delta}\,\sqrt{1+\log_2 (s^{-1})}$ is increasing on $(0,1]$; thus taking $c=1$, it follows that the lemma holds. When $\overline{s}_\delta\in (0,1)$, then the latter function is increasing on $(0,\overline{s}_\delta]$ and decreasing on $[\overline{s}_\delta\,1]$; thus setting
$$
c=\frac{\max_{s\in (0,1]}\big\{s^{\delta}\,\sqrt{1+\log_2 (s^{-1})}\big\}}
{\min_{s\in [\overline{s}_\delta,1]}\big\{s^{\delta}\,\sqrt{1+\log_2 (s^{-1})}\big\}},
$$
one obtains the lemma. $\square$\\

Now we are in position to prove Theorem~\ref{thm:z}.\\
\noindent{\bf Proof of Theorem~\ref{thm:z}:} First notice that it is sufficient to show that
for all $\omega\in\Omega^*$ (recall that $\Omega^*$ is the event of probability~$1$ on which (\ref{lambda}) holds), $J\in\N$ and $Q\in\N\setminus\{0\}$ one has,
\begin{equation}
\label{eq30:ant-ch4}
\big \|Z_{J+Q}(\cdot,\omega)-Z_J (\cdot,\omega)\big\|_{C^\gamma ([0,1])}\le D(\omega) 2^{-J\min(1/2-\gamma,\alpha-1/2)}\,\sqrt{1+J}.
\end{equation}
Indeed, (\ref{eq30:ant-ch4}) implies that $\big(Z_J (\cdot,\omega)\big)_{J\in\N}$ is a Cauchy
sequence in the Banach space $C^\gamma ([0,1])$ and, as consequence, it converges, in this space, to some limit
denoted by $\widetilde{Z}(\cdot,\omega)$. Next let $\Omega_{2}^*$ be the event
of probability $1$ defined as $$\Omega_{2}^*=\bigcap_{q\in [0,1]\cap \Q} \Omega_{1,q}^*,$$ where
$\Omega_{1,q}^*\subseteq \Omega^*$ is the event $\Omega_{1,t}^*$ introduced in Lemma~\ref{lem1:ant-ch4} when $t=q$. Thus it follows from the latter lemma, that for each $\omega\in\Omega_{2}^*$ and all $q\in[0,1]\cap \Q$ one has $Z(q,\omega)=\widetilde{Z}(q,\omega)$; this is equivalent to $Z(\cdot,\omega)=\widetilde{Z}(\cdot,\omega)$, since $Z(\cdot,\omega)$ and $\widetilde{Z}(\cdot,\omega)$ are continuous function. Then, letting in (\ref{eq30:ant-ch4}), $Q$ goes to $+\infty$ while $J$ fixed, one gets (\ref{eq20:ant-ch4}).

From now on, our goal will be to show that (\ref{eq30:ant-ch4}) holds. To this end, in view of (\ref{eq:holdnorm-ch4}), it is sufficient to prove that there exist two positive random variables $D_1$ and $D_2$ of finite moment of any order, such that one has for all $\omega\in\Omega^*$, $J\in\N$ and
$Q\in \N\setminus\{0\}$,
\begin{equation}
\label{eq31:ant-ch4}
\begin{split}
& \sup_{0\le t_1 <t_2\le 1} \frac{\big |Z_{J+Q}(t_1,\omega)-Z_J (t_1,\omega)-Z_{J+Q}(t_2,\omega)+Z_J (t_2,\omega)\big |}{|t_1-t_2|^\gamma}\\
& \le D_1(\omega) 2^{-J\min(1/2-\gamma,\alpha-1/2)}\,\sqrt{1+J}
\end{split}
\end{equation}
and
\begin{equation}
\label{eq32:ant-ch4}
\sup_{0\le t\le 1} \big |Z_{J+Q}(t,\omega)-Z_J (t,\omega)\big|\le D_2(\omega) 2^{-J\min(1/2-\gamma,\alpha-1/2)}\,\sqrt{1+J}.
\end{equation}
We will only show that (\ref{eq31:ant-ch4}) is satisfied, since (\ref{eq32:ant-ch4}) can be
obtained in the same way. Let $t_1<t_2$ be two arbitrary and fixed real numbers belonging
to the interval $[0,1]$. Observe that in view of (\ref{eq15:ant-ch4}), using the triangle inequality, one has,
\begin{eqnarray}
\label{eq33:ant-ch4}
\nonumber
&& |Z_{J+Q}(t_1,\omega)-Z_J (t_1,\omega)-Z_{J+Q}(t_2,\omega)+Z_J (t_2,\omega)\big |\\
\nonumber
&& =\Big |\sum_{j=J+1}^{J+Q}\sum_{k=0}^{2^j -1}\big (a_{j,k}(t_1,\omega)-a_{j,k}(t_2,\omega)
\big )\lambda_{j,k}(\omega)\Big |\\
&& \le \sum_{j=J+1}^{+\infty}\Big |\sum_{k=0}^{2^j -1}\big (a_{j,k}(t_1,\omega)-a_{j,k}(t_2,\omega)
\big )\lambda_{j,k}(\omega)\Big |.
\end{eqnarray}
Let us now give appropriate bounds for
\begin{equation}
\label{eq34:ant-ch4}
U_j (t_1,t_2,\omega)=\Big |\sum_{k=0}^{2^j -1}\big (a_{j,k}(t_1,\omega)-a_{j,k}(t_2,\omega)
\big )\lambda_{j,k}(\omega)\Big |
\end{equation}
and derive form them (\ref{eq31:ant-ch4}). We denote by $j_0\in\N$ the unique integer such that
\begin{equation}
\label{eq35:ant-ch4}
2^{-j_0-1}<|t_1-t_2|\le 2^{-j_0}.
\end{equation}
Also for all $t\in [0,1]$ and $j\in\N$, we denote by $\tilde{k}(j,t)$ the unique integer
belonging to $\{0,\ldots, 2^j-1\}$ such that
\begin{equation}
\label{eq38:ant-ch4}
\frac{\tilde{k}(j,t)}{2^j}\le t <\frac{\tilde{k}(j,t)+1}{2^j},
\end{equation}
with the convention that
\begin{equation}
\label{eq39:ant-ch4}
\tilde{k}(j,1)=2^j-1.
\end{equation}
Observe that when $t\in [0,1)$,
\begin{equation}
\label{eq40:ant-ch4}
\tilde{k}(j,t)=\big [2^j t\big],
\end{equation}
$[\cdot]$ being the integer part function. Also, observe that (\ref{eq34:ant-ch4}), (\ref{eq25:ant-ch4}), and the fact that $a_{j,k}(t,\omega)$ does not $t$ when $t\ge (k+1)/2^j$ (see (\ref{eq27:ant-ch4})), imply that
\begin{equation}
\label{eq36:ant-ch4}
U_j (t_1,t_2,\omega)=\Big |\sum_{k=\tilde{k}(j,t_1)}^{\tilde{k}(j,t_2)} \big (a_{j,k}(t_1,\omega)-a_{j,k}(t_2,\omega)
\big )\lambda_{j,k}(\omega)\Big |.
\end{equation}
Let us now study the following two cases: $j\le j_0$ and $j>j_0$. First assume that $j\le j_0$, then
(\ref{eq35:ant-ch4}) implies that $2^{-j}\ge |t_1-t_2|$ and, as a consequence that
$\tilde{k}(j,t_2)\in\big\{\tilde{k}(j,t_1),\tilde{k}(j,t_1)+1\big\}$. When,
$\tilde{k}(j,t_2)=\tilde{k}(j,t_1)$, it follows from (\ref{eq36:ant-ch4}), (\ref{lambda}) and (\ref{eq19:ant-ch4}), that
\begin{eqnarray*}
U_j (t_1,t_2,\omega)&=&\big |a_{j,\tilde{k}(j,t_1)}(t_1,\omega)-a_{j,\tilde{k}(j,t_1)}(t_2,\omega)
\big |\big|\lambda_{j,\tilde{k}(j,t_1)}(\omega)\big|\\
&\le & G_1(\omega) |t_1-t_2|2^{j/2}\,\sqrt{1+j},
\end{eqnarray*}
where $G_1(\omega):=C_{*}(\omega)C_1(\omega)$. When, $\tilde{k}(j,t_2)=\tilde{k}(j,t_1)+1$, putting together (\ref{eq36:ant-ch4}), (\ref{lambda}), (\ref{eq19:ant-ch4}), (\ref{eq24:ant-ch4}) and the fact that $(\tilde{k}(j,t_1)+1)/2^j \in [t_1,t_2]$, one obtains that,
\begin{equation*}
\begin{split}
& U_j (t_1,t_2,\omega)\\
& \le \big |a_{j,\tilde{k}(j,t_1)}(t_1,\omega)-a_{j,\tilde{k}(j,t_1)}(t_2,\omega)
\big |\big|\lambda_{j,\tilde{k}(j,t_1)}(\omega)\big|+\big|a_{j,\tilde{k}(j,t_1)+1}(t_2,\omega)
\big |\big|\lambda_{j,\tilde{k}(j,t_1)+1}(\omega)\big|\\
& \le G_1 (\omega) |t_1-t_2|2^{j/2}\,\sqrt{1+j}+G_2(\omega) \big(t_2-(\tilde{k}(j,t_1)+1)/2^j\big)|t_1-t_2|2^{j/2}\,\sqrt{1+j}\\
& \le G_3(\omega)|t_1-t_2|2^{j/2}\,\sqrt{1+j},
\end{split}
\end{equation*}
where $G_2(\omega):= A(\omega)C_{*}(\omega)$ and $G_3(\omega):=G_1(\omega)+G_2(\omega)$. Thus, we have
shown that there is a positive random variable $G_4$ of finite moment of any order, non depending on $\gamma$, $j_0$, $t_1$ and $t_2$, such that for all $j\le
j_0$, one has,
$$
\frac{U_j (t_1,t_2,\omega)}{|t_1-t_2|^\gamma}\le G_4(\omega)|t_1-t_2|^{1-\gamma}\, 2^{j/2}\,\sqrt{1+j}.
$$
Therefore, in view of (\ref{eq35:ant-ch4}), for each integer $J$ satisfying $0\le J < j_0$, one has,
\begin{equation}
\label{eq37:ant-ch4}
\begin{split}
\sum_{j=J+1}^{j_0}\frac{U_j (t_1,t_2,\omega)}{|t_1-t_2|^\gamma} &\le G_4(\omega)|t_1-t_2|^{1-\gamma}\, \sum_{j=J+1}^{j_0} 2^{j/2}\,\sqrt{1+j}\\
&\le  \sqrt{2}\big(\sqrt{2}-1\big)^{-1}G_4(\omega)|t_1-t_2|^{1-\gamma}\, 2^{j_0/2}\,\sqrt{1+j_0}\\
&\le \sqrt{2}\big(\sqrt{2}-1\big)^{-1}G_4(\omega)|t_1-t_2|^{1/2-\gamma}\,\sqrt{1+\log_2 (|t_1-t_2|^{-1})}.
\end{split}
\end{equation}
Next, it follows from (\ref{eq37:ant-ch4}), the inequality $2^{-J}>2^{-j_0}\ge |t_1-t_2|$ (see (\ref{eq35:ant-ch4})) and Lemma~\ref{lem3:ant-ch4} (in which one
takes $\delta=1/2-\gamma$, $x=|t_1-t_2|$ and $y=2^{-J}$), that
\begin{equation}
\label{eq44:ant-ch4}
\sum_{j=J+1}^{j_0}\frac{U_j (t_1,t_2,\omega)}{|t_1-t_2|^\gamma}\le G_5(\omega) 2^{-J(1/2-\gamma )}\sqrt{1+J},
\end{equation}
where
$$
G_5(\omega):=c\sqrt{2}\big(\sqrt{2}-1\big)^{-1}G_4(\omega),
$$
$c$ being the constant introduced in Lemma~\ref{lem3:ant-ch4}. Let us now study the case where $j>j_0$. Observe that in this case, in view of Relations (\ref{eq35:ant-ch4}) and (\ref{eq38:ant-ch4}),  one necessarily has $\tilde{k}(j,t_1)<\tilde{k}(j,t_2)$ and
\begin{equation}
\label{eq41:ant-ch4}
2^{j}|t_1-t_2|> 1.
\end{equation}
Also observe that, using (\ref{eq39:ant-ch4}), (\ref{eq40:ant-ch4}) and (\ref{eq41:ant-ch4}), one obtain that
\begin{equation}
\label{eq42:ant-ch4}
\tilde{k}(j,t_2)-\tilde{k}(j,t_1)< 2^j t_2 -2^j t_1 +1<2^{j+1} |t_1-t_2|.
\end{equation}
It follows from (\ref{eq34:ant-ch4}), (\ref{eq27:ant-ch4}), (\ref{eq38:ant-ch4}), Lemma~\ref{lem2:ant-ch4}, (\ref{eq42:ant-ch4}), (\ref{lambda}) and (\ref{eq19:ant-ch4}), that
\begin{equation}
\label{eq43:ant-ch4}
\begin{split}
U_j (t_1,t_2,\omega) & \le \big |a_{j,\tilde{k}(j,t_1)}(t_1,\omega)\big |\big|\lambda_{j,\tilde{k}(j,t_1)}(\omega)\big|+\sum_{k=\tilde{k}(j,t_1)}^{\tilde{k}(j,t_2)-1}
\big |a_{j,k}(t_2,\omega)|\big|\lambda_{j,k}(\omega)\big|\\
&\quad+\big |a_{j,\tilde{k}(j,t_2)}(t_2,\omega)\big|\big|\lambda_{j,\tilde{k}(j,t_2)}(\omega)\big|\\
& \le  2G_2 (\omega)\Big (2^{-j/2}\,\sqrt{1+j} + |t_1-t_2| \,2^{-j(\alpha-1/2)}\,\sqrt{1+j}\Big).
\end{split}
\end{equation}
Thus (\ref{eq43:ant-ch4}) and (\ref{eq41:ant-ch4}) imply that for all $j> j_0$, one has,
\begin{eqnarray}
\label{eq45:ant-ch4}
\nonumber
\frac{U_j (t_1,t_2,\omega)}{|t_1-t_2|^\gamma} & \le & 2G_2 (\omega)\Big (2^{-j(1/2-\gamma)}\,\sqrt{1+j} +|t_1-t_2|^{1-\gamma} \, 2^{-j(\alpha-1/2)}\,\sqrt{1+j}\Big)\\
&\le & G_6 (\omega) 2^{-j\min(1/2-\gamma,\alpha-1/2)}\,\sqrt{1+j}.
\end{eqnarray}
where $G_6 (\omega):=4G_2(\omega)$. Next, observe that one has for $J\in\N$,
\begin{eqnarray}
\label{eq46:ant-ch4}
\nonumber
&&\sum_{j=J+1}^{+\infty}2^{-j\min(1/2-\gamma,\alpha-1/2)}\,\sqrt{1+j}\\
\nonumber
&&=2^{-J\min(1/2-\gamma,\alpha-1/2)}\,\sqrt{1+J}\sum_{j=J+1}^{+\infty}
2^{-(j-J)\min(1/2-\gamma,\alpha-1/2)}\,\sqrt{\frac{1+j}{1+J}}\\
&& \le c_1 2^{-J\min(1/2-\gamma,\alpha-1/2)}\,\sqrt{1+J},
\end{eqnarray}
where the constant
$$
c_1:=\sum_{l=1}^{+\infty} 2^{-l\min(1/2-\gamma,\alpha-1/2)}\,\sqrt{1+l}<+\infty.
$$
Thus combining (\ref{eq45:ant-ch4}) with (\ref{eq46:ant-ch4}), it follows that for all $J\ge j_0$, one has,
\begin{equation}
\label{eq47:ant-ch4}
\sum_{j=J+1}^{+\infty} \frac{U_j (t_1,t_2,\omega)}{|t_1-t_2|^\gamma} \le G_7 (\omega) 2^{-J\min(1/2-\gamma,\alpha-1/2)}\,\sqrt{1+J},
\end{equation}
where $G_7 (\omega):=c_1G_6(\omega)$. Let us now show that, for all $J\in\N$, one has,
\begin{equation}
\label{eq48:ant-ch4}
\sum_{j=J+1}^{+\infty} \frac{U_j (t_1,t_2,\omega)}{|t_1-t_2|^\gamma} \le D_1 (\omega) 2^{-J\min(1/2-\gamma,\alpha-1/2)}\,\sqrt{1+J},
\end{equation}
where $D_1(\omega):=(1+c)G_7 (\omega)+G_5(\omega)$, $c$ is the constant introduced in Lemma~\ref{lem3:ant-ch4} and $G_5(\omega)$ has been introduced
in (\ref{eq44:ant-ch4}). It is clear that (\ref{eq47:ant-ch4}) implies
that (\ref{eq48:ant-ch4}) holds when $J\ge j_0$, so from now on, we assume that $j_0\ge 1$ and that $J$ is an arbitrary nonnegative integer satisfying $J<j_0$.
It follows from (\ref{eq47:ant-ch4}) that
$$
\sum_{j=j_0+1}^{+\infty} \frac{U_j (t_1,t_2,\omega)}{|t_1-t_2|^\gamma} \le G_7 (\omega) 2^{-j_0\min(1/2-\gamma,\alpha-1/2)}\,\sqrt{1+j_0}.
$$
Then using Lemma~\ref{lem3:ant-ch4} (in which one takes $\delta=\min(1/2-\gamma,\alpha-1/2)$, $x=2^{-j_0}$ and $y=2^{-J}$), one obtains that
\begin{equation}
\label{eq49:ant-ch4}
\sum_{j=j_0+1}^{+\infty} \frac{U_j (t_1,t_2,\omega)}{|t_1-t_2|^\gamma} \le c G_7 (\omega)2^{-J\min(1/2-\gamma,\alpha-1/2)}\,\sqrt{1+J}.
\end{equation}
Next combining (\ref{eq49:ant-ch4}) with (\ref{eq44:ant-ch4}), it follows that (\ref{eq48:ant-ch4}) holds in the case where $J<j_0$. Finally
(\ref{eq33:ant-ch4}), (\ref{eq34:ant-ch4}) and (\ref{eq48:ant-ch4}) imply that (\ref{eq31:ant-ch4}) is satisfied. $\square$

\section{Simulation of $\{Z(t)\}_{t\in[0,1]}$ via the Haar multiresolution analysis}
\label{sec:ch4-simulation}
Our algorithm for simulating $\{Z(t)\}_{t\in[0,1]}$ mainly relies on Theorem~\ref{thm:z} which allows to approximate $\{Z(t)\}_{t\in[0,1]}$, for $J$ large enough,
by the process $\{Z_J(t)\}_{t\in[0,1]}$ defined in (\ref{eq15:ant-ch4}) as a finite of sum. First we will give an expression for the latter process which makes it rather easy to simulate, to this end, we need to introduce some notations.
\begin{itemize}
\item For all fixed integers $J\in\N$ and $l\in\{0,\ldots , 2^J -1\}$, the function $\varphi_{J,l}$ is defined as,
\begin{equation}
\label{eq50:ant-ch4}
\varphi_{J,l}=2^{J/2}\mathds{1}_{[\frac{l}{2^J},\frac{l+1}{2^J})}.
\end{equation}
\item For all fixed integer $J\in\N$, we denote by $\big\{\delta_{J,l}:\,l\in\{0,\ldots, 2^{J}-1\}\big\}$ the finite sequence of independent standard
Gaussian random variables defined as,
\begin{equation}
\label{eq51:ant-ch4}
\delta_{J,l}:=\int_0^1\varphi_{J,l}(s)\ud W(s)=2^{J/2}\Big (W(\frac{l+1}{2^J})-W(\frac{l}{2^J})\Big ).
\end{equation}
\item For all fixed integers $J\in \N$ and $l\in\{0,\dots, 2^J-1\}$, the stochastic process $\{b_{J,l}(t)\}_{t\in [0,1]}$ is defined for all
$(t,\omega)\in [0,1]\times\Omega$, as,
\begin{equation}
\label{eq51bis:ant-ch4}
b_{J,l}(t,\omega)=\int_0^1 K(t,s,\omega)\varphi_{J,l}(s)\ud s,
\end{equation}
where $K(t,s,\omega)$ has been introduced in (\ref{eq5:ant-ch4}).
\end{itemize}
The following proposition provides a nice expression for the process $\{Z_J (t)\}_{t\in [0,1]}$.
\begin{prop}
\label{prop1:ant-ch4}
For all fixed integer $J\in\N$, let $\{Z_J (t)\}_{t\in [0,1]}$ be the stochastic process defined in (\ref{eq15:ant-ch4}). Then one has for all $t\in [0,1]$, almost
surely,
\begin{equation}
\label{eq52:ant-ch4}
Z_J(t)=\sum_{l=0}^{2^{J+1}-1} b_{J+1,l}(t)\delta_{J+1,l}.
\end{equation}
\end{prop}

Observe that, Relation (\ref{eq52:ant-ch4}), also holds almost surely for all $J\in\N$ and $t\in [0,1]$ (i.e. this relation is satisfied on an event of 
probability~$1$ which does not depend on $J$ and $t$), since the trajectories of the processes $\{Z_J(t)\}_{t\in [0,1]}$
and $\left\{\sum_{l=0}^{2^{J+1}-1} b_{J+1,l}(t)\delta_{J+1,l}\right\}_{t\in [0,1]}$, are with probability $1$, continuous functions.\\

\noindent {\bf Proof of Proposition~\ref{prop1:ant-ch4}:} The main ingredient of the proof is the Haar multiresolution analysis (see for example \cite{M1,D,W}) of the Hilbert space $L^2([0,1])$; that is the increasing (in the sense of the inclusion) sequence $(V_J)_{J\in\Z_+}$ of the finite dimensional subspaces of $L^2([0,1])$ defined as, 
$$
V_0:=\mbox{Span}\,\big\{\varphi_{0,0}\big\}
$$
 and for all integer $J\ge 1$,
$$
V_J:=\mbox{Span}\,\big\{\varphi_{0,0},\,\psi_{j,k}:\,\,j\in\{0,\ldots, J-1\}\mbox{ and } k\in \{0,\ldots, 2^j-1\}\big\},
$$
where the orthonormal functions $\varphi_{0,0}$ and $\psi_{j,k}$ have been introduced respectively in (\ref{eq16:ant-ch4}) and in (\ref{eq17:ant-ch4}).
Relations (\ref{dec:K}) and (\ref{eq9bis:ant-ch4}) imply that for every fixed $J\in\Z_+$ and $(t,\omega)\in [0,1]\times \Omega$, the function
$K_J(t,\cdot,\omega)$, can be viewed as the orthogonal projection (in the sense of the usual inner product of $L^2([0,1])$) of the function $K(t,\cdot,\omega)$, on the space $V_{J+1}$. On the other hand, it is known (see for example \cite{M1,D,W}) that, the finite sequence, $$\big\{\varphi_{J+1,l},\,l\in \{0,\ldots, 2^{J+1}-1\}\big\},$$
forms an orthonormal basis of $V_{J+1}$. Therefore, in view of (\ref{eq51bis:ant-ch4}), one has, for all $s\in [0,1]$,
\begin{equation}
\label{eq54:ant-ch4}
K_J(t,s,\omega)=\sum_{l=0}^{2^{J+1} -1} b_{J+1,l}(t,\omega)\varphi_{J+1,l}(s).
\end{equation}
Finally, using the fact that,
$$
Z_J (t)=\int_0^1K_J(t,s)\ud W(s)
$$
as well as Relations (\ref{eq9bis:ant-ch4}), (\ref{eq14:ant-ch4}), (\ref{eq54:ant-ch4}) and (\ref{eq51:ant-ch4}),
one obtains (\ref{eq52:ant-ch4}). $\square$\\

Now, we are ready to describe the main step of our algorithm for simulating $\{Z(t)\}_{t\in [0,1]}$.\\
\\
{\bf Main steps of our algorithm for simulating $\{Z(t)\}_{t\in [0,1]}$:}

\begin{itemize}
\item[(1)] We take $J$ large enough and we simulate the finite sequence
$$\big\{\delta_{J,l}:\,l\in\{0,\ldots, 2^{J}-1\}\big\}$$
of the standard independent Gaussian random variables defined in~(\ref{eq51:ant-ch4}).
\item[(2)] We simulate
$$\left\{X(0),X\Big(\frac{1}{2^{2J}}\Big),\ldots, X\Big(\frac{2^{2J}-1}{2^{2J}}\Big)\right\};$$
observe that in the case where the centered Gaussian process $\{X(t)\}_{t\in [0,1]}$ is a multifractional Brownian motion; such a simulation can be done by making use of one of the efficient methods described in \cite{B1}.
\item[(3)] Noticing that, for all $l\in\{0,\ldots, 2^J-1\}$, and $m\in\{l+1,\ldots, 2^J\}$, Relations (\ref{eq50:ant-ch4}), (\ref{eq51bis:ant-ch4}) and (\ref{eq5:ant-ch4}) imply that,
$$
b_{J,l}\Big(\frac{m}{2^J}\Big)=2^{J/2} \int_{\frac{l}{2^J}}^\frac{l+1}{2^J}\Phi(X(s))\ud s,
$$
we approximate the latter integral, by the Riemann sum
\begin{equation}
\label{eq55:ant-ch4}
\widehat{b}_{J,l}\Big(\frac{m}{2^J}\Big):=2^{-3J/2}\sum_{q=0}^{2^J-1}\Phi\left(X\Big(\frac{l}{2^J}+\frac{q}{2^{2J}}\Big)\right).
\end{equation}
On the other hand, observe that (\ref{eq50:ant-ch4}), (\ref{eq51bis:ant-ch4}) and (\ref{eq5:ant-ch4}) entail that for each $m\in \{0,\ldots, l\}$,
$$
b_{J,l}\Big(\frac{m}{2^J}\Big)=0
$$
\item[(4)] Thus, in view of (\ref{eq52:ant-ch4}), for all $m\in\{1,\ldots, 2^J\}$, we approximate $Z_{J-1}\big (\frac{m}{2^J}\big)$ by
\begin{equation}
\label{eq55bis:ant-ch4}
\widehat{Z}_{J-1}\Big (\frac{m}{2^J}\Big):=\sum_{l=0}^{m-1} \widehat{b}_{J,l}\Big(\frac{m}{2^J}\Big) \delta_{J,l}.
\end{equation}
Then we simulate
$$
\left\{\widehat{Z}_{J-1}\Big (\frac{1}{2^J}\Big),\ldots, \widehat{Z}_{J-1}\Big (\frac{2^J -1}{2^J}\Big), \widehat{Z}_J(1)\right\},
$$
by using the fact that
\begin{equation}
\label{eq56:ant-ch4}
\widehat{Z}_{J-1}\Big (\frac{1}{2^J}\Big)=\widehat{b}_{J,0}\Big(\frac{1}{2^J}\Big) \delta_{J,0}
\end{equation}
and the induction relation, for all $m\in\{2,\ldots, 2^J\}$,
\begin{equation}
\label{eq57:ant-ch4}
\widehat{Z}_{J-1}\Big (\frac{m}{2^J}\Big)=\widehat{Z}_{J-1}\Big (\frac{m-1}{2^J}\Big)+\widehat{b}_{J,m-1}\Big(\frac{m}{2^J}\Big) \delta_{J,m-1}.
\end{equation}
Observe that (\ref{eq56:ant-ch4}) and (\ref{eq57:ant-ch4}) easily result from (\ref{eq55:ant-ch4}) and (\ref{eq55bis:ant-ch4}).
\item[(5)] Finally, by interpolating the $2^J +1$ points
$$
(0,0)\,;\, \left(\frac{1}{2^J}, \widehat{Z}_{J-1}\Big (\frac{1}{2^J}\Big)\right)\,;\,\ldots\,;\, \left(\frac{2^J-1}{2^J},\widehat{Z}_{J-1}\Big (\frac{2^J-1}{2^J}\Big)\right)\,;\,
\big(1,\widehat{Z}_{J-1}(1)\big)
$$
we obtain a stochastic process $\big\{\widehat{Z}_{J-1}(t)\big\}_{t\in [0,1]}$ which, in view of Theorem~\ref{thm:z}, satisfies the following property: for all fixed
$\gamma\in [0,1/2)$, there exists a random variable $D'$ of finite moment of any order, such that one has, almost surely, for all $J$ big enough,
$$
\big \|Z-\widehat{Z}_{J-1} \big\|_{C^\gamma ([0,1])}\le D' 2^{-J\min(1/2-\gamma,\alpha-1/2)}\,\sqrt{1+J},
$$
where $\|\cdot\|_{C^\gamma ([0,1])}$ is the usual norm on the H\"older space $C^\gamma ([0,1])$.
\end{itemize}
$\square$

\newpage 

\begin{center}
\begin{figure}
\includegraphics[width=\textwidth]{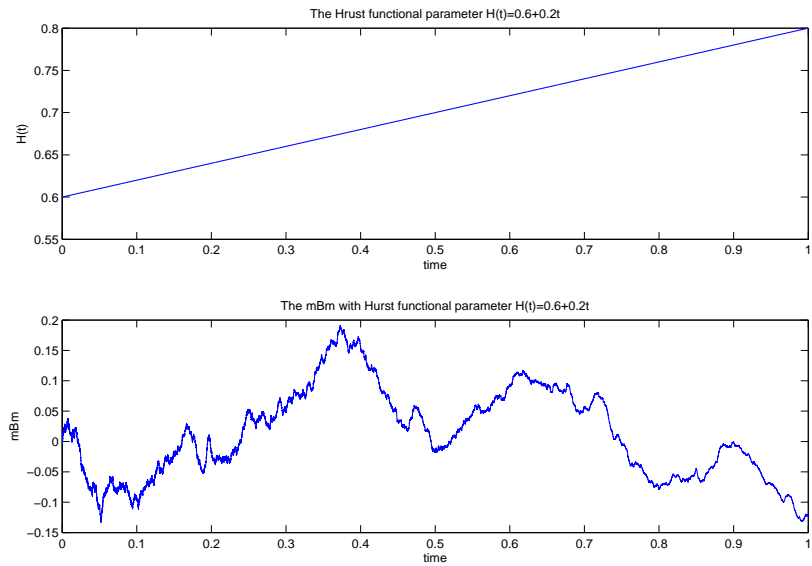}
\caption{Graph of the function $H(\cdot)$ such that $H(s)=0.6+0.2s$ for all $s\in [0,1]$, and a simulation of a trajectory of a multifractional Brownian motion $\{X(s)\}_{s\in [0,1]}$ with Hurst functional parameter $H(\cdot)$.}
\label{image1}
\end{figure}
\end{center}
\begin{center}
\begin{figure}
\includegraphics[width=\textwidth]{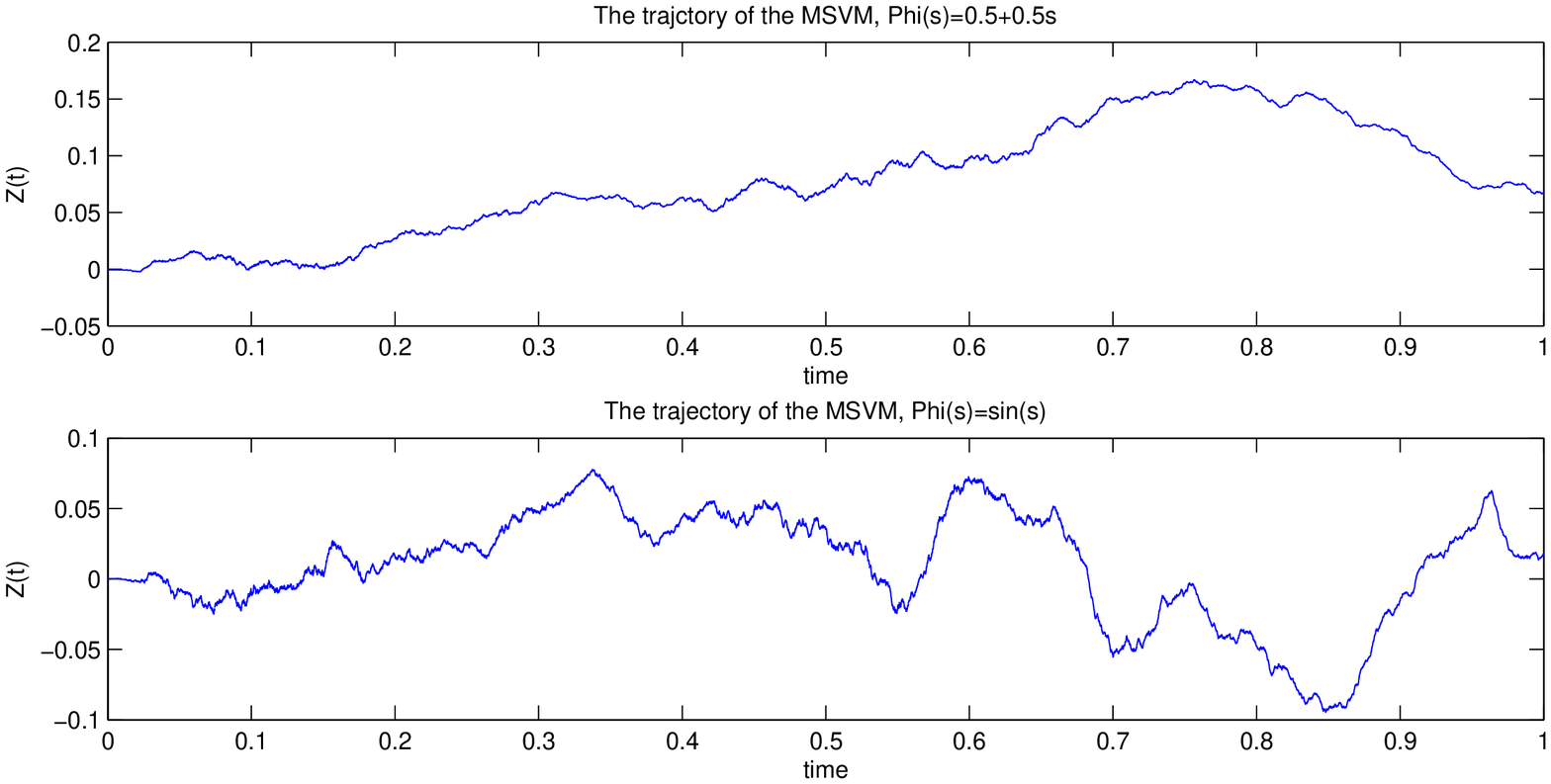}
\caption{Simulation, in the case where $\Phi(x)=0.5+0.5x$ for all $x\in\R$, and also in the case where 
$\Phi(x)=\sin(x)$ for all $x\in\R$, of a trajectory of the log price process $\{Z(t)\}_{t\in [0,1]}$ generated through (\ref{eq0:ant-ch4}) by the multifractional Brownian motion $\{X(s)\}_{s\in [0,1]}$
in Figure~\ref{image1}. The shapes of these two simulated trajectories tend to confirm the fact that the pointwise H\"older exponent of $\{Z(t)\}_{t\in (0,1]}$ does not change from one place to another and is equal to $1/2$ (see Theorem~\ref{thm:1}).}
\end{figure}
\end{center}

\newpage 

\bibliographystyle{plain}
\bibliography{Thesisbib}

\begin{quote}
\begin{small}
\textsc{Antoine Ayache and Qidi Peng}: U.M.R. CNRS 8524, Laboratoire Paul
Painlev\'e, B\^at. M2, Universit\'e Lille 1, 59655 Villeneuve d'Ascq Cedex, France.\\
E-mails: \texttt{Antoine.Ayache@math.univ-lille1.fr, \ Qidi.Peng@math.univ-lille1.fr}
\end{small}
\end{quote}

\end{document}